# On a Fundamental Task of Bidiophantine Geometric Figures

Zurab Agdgomelashvili

(Georgian Technical University, 77, M. Kostava str.,

Tbilisi, 0175, Georgia)

**Abstract:** The goal of the work is to take on and study one of the fundamental tasks studying Bidiophantine polygons (let us call a polygon Diophantine, if the distance between each two vertex of those is expressed by a natural number and we say that a Diophantine polygon is Bidiophantine if the coordinates of its each vertex are integer numbers).

Task **($n; k$): is there a Bidiophantine n-gon ($n \geq 3$) with any side or diagonal equal to $k$ for each fixed natural number $k$ ? In case if it exists then let us find each such $n$.

As a result of fundamental studies we have obtained a full answer to the above mentioned question for k=1, k=2, k=3 and k=4.

The paper shows that for k=1 and k=2 such Bidiophantine n-gon does not exist, and for k=3 as well, as for k=4 definitely n is equal to 3 or 4. It is also shown that when $k>2$ there is always a Bidiophantine rectangle with the length of one of the sides equal to k.

## INTRODUCTION

Let us call a polygon Diophantine, if the distance between each two vertex of those is expressed by a natural number.

Let us call as bidiophantine a Diophantine polygon those coordinates of each vertex are integers.

In this paper, we will consider one of assigned by us task that undoubtedly represents one of the fundamental tasks of studying the properties of bidiophantine polygons.

**Task** **($n;k$)** exists if for not each fixed k natural number, bidiophantine n-gon ($n≥3$) whose length of arbitrary side or diagonal is equal to k, and if it exists, then let's find all such *n*.

Let's firstly consider the task* (n; k).

From the tasks of studying Diophantine geometric figures, one of the most important obviously is the following:

**Task* (n;k) exists or not for each fixed k natural number, Diophantine n-gon ($n≥3$) whose length of arbitrary side or diagonal is equal to k, and if it exists, then let's find all such *n*.**

It is shown by us that there does not exist such Diophantine n-gon, both for convex and concave, whose length of arbitrary sides or diagonals are equal to 1. I.e. the above issue has been solved for $k = 1$.

**BASIC PART**

To prove this, let's consider a few tasks.

**Lemma 1. If lengths of one of the side of Diophantine triangle is equal to 1, then the rest of its sides represent equal legs of the triangle.**

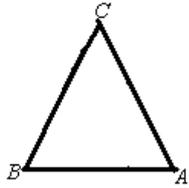
Fig. 1

Task: $\triangle ABC$;

$|BC|, |AC| \in N$;

$|AB| = 1$.

u.d : $|BC| = |AC|$.

Without limiting of generality let us assume that $|BC| \leq |AC|$. By the inequality of triangle $|AB| + |BC| > |AC|$ or $1 + |BC| > |AC|$.

I.e. $\begin{cases} |BC|, |AC| \in N; \\ |BC| \leq |AC| < |BC| + 1. \end{cases} \Rightarrow (|BC| = |AC|)$. Q.E.D.

This remarkable property is one of the cornerstones of the Diophantine and bidiophantine geometric figure's research apparatus.

**Lemma 2. The length of each side and each diagonal of the convex Diophantine rectangle is greater than 1.**

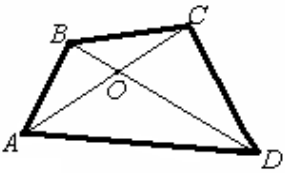
Fig.2

Let's assume the opposite. I.e. let's say that in rectangle *ABCD* the length of each side and diagonal is expressed by the natural number, and also the length of one of the sides, without limiting of generality let's say that $|AB| = 1$.

**According to Lemma 1** $|BC| = |AC|$ and $|BD| = |AD|$.

This is impossible, since then *C* and *D* the points should be located on the median [*AB*]. I.e. our assumption is false. Therefore acceleration of this condition, the length of each side of $\square ABCD$ is greater than 1.

Now let's say that length of the $\square ABCD$ diagonal is equal to 1, without limitation of generality let's say $|AC| = 1$.

According to Lemma 1 $|AB| = |BC|$ and $|CD| = |AD|$. Using the inequality of triangle $\triangle BOC$ and $\triangle AOD$, it is easy to show, that $|BD| + |AC| > |BC| + |AD|$ or $|AB| + |AD| < |BD| + 1$. From $\triangle ABD$ we have $|AB| + |AD| > |BD|$. I.e. $\begin{cases} |AB|, |BD|, |AB| \in N; \\ |BD| < |AB| + |AD| < |BD| + 1. \end{cases}$ is impossible. Thus the assumption is false, or $|AC| > 1$.

**I.e. finally we have that the length of each side and diagonal of each convex Diophantine rectangle is greater than 1. Q.E.D.**

**Lemma 3. If none of the located on the plane four points are not located on one straight line and at the same time distance between each of them is expressed by a natural number, then each of these distances is greater than 1.**

Let's assume the opposite. I.e. let's say that there are four points on the plane, none of that are located on one straight line, the distance between each of them is expressed by a natural number, and at the same time some of them are equal to 1. Then, according to Lemma 1, the remaining both points must be located on a median of segment with equal to 1 length. Here we would have two cases

1) 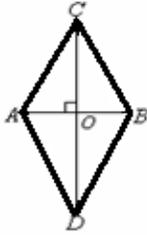 Task: $\triangle ABC$; $|AB|=1$; $|AC|=|BC| \in N$; $|AD|=|DB| \in N$; $|CD| \in N$.

Fig.3

2) 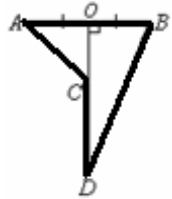 Task: $|AB|=1$, $|AC|=|BC|=m$; $|AD|=|DB|=n$, $|CD|=l$, $m,n,l \in N$.

Fig.4

1) This case is considered in the **Lemma 2**.

2) From $\triangle AOC$ and $\triangle DOB$ we have:

$$\begin{cases} |AC|^2 = |AO|^2 + |OC|^2; \\ |BD|^2 = |OB|^2 + |OD|^2; \\ |OD| = |OC| + |CD|, |AO| = |OB| = 0{,}5; \\ |AC| = |CB| = m, |AD| = |DB| = n, |CD| = l, m,n,l \in N. \end{cases} \Rightarrow$$

$$\Rightarrow \begin{cases} m^2 = \left(\frac{1}{2}\right)^2 + |OC|^2; \\ n^2 = \left(\frac{1}{2}\right)^2 + (l + |OC|)^2; \\ m,n,l \in N. \end{cases} \Rightarrow \begin{cases} |OC| = \dfrac{n^2 - m^2 - l^2}{2l}; \\ 4m^2 - 1 = (2|OC|)^2; \\ m,n,l \in N. \end{cases} \Rightarrow \begin{cases} |OC| = \dfrac{n^2 - m^2 - l^2}{l} = q \in Q_+; \\ 4m^2 - 1 = q^2; \\ m,n,l \in N. \end{cases} \Rightarrow$$

$$\Rightarrow \begin{cases} m,q \in N; \\ 4m^2 - q^2 = 1. \end{cases} \Rightarrow \begin{cases} m,q \in N; \\ (2m-q)(2m+q) = 1. \end{cases} \Rightarrow \begin{cases} m,q \in N; \\ 2m - q = 1; \\ 2m + q = 1. \end{cases}$$

This is impossible. Thus the distance between arbitrary two points from these given points is greater than 1.

**Theorem 1. The length of each side and each diagonal of convex Diophantine *n*-gon ($n > 3$) is greater than 1.**

Let's assume the opposite. I.e. assume a convex Diophantine *n*-polygon ($n > 3$), the distance between any two vertexes of that is equal to 1. By virtue of Lemma 1, the remaining vertexes should be located on the median of the connecting these two vertexes segment. Because this *n* polygon should be convex, so obviously is *n*= 4, but according to Lemma 2 such rectangle does not exist. I.e. our assumption is false, or the length of each side and each diagonal of arbitrary convex Diophantine *n*-gon ($n > 3$) is greater than 1. Q.E.D.

**Theorem 2. If from located on the plane n-point ($n > 3$), none of the three point is located on one straight line and the distance between each of them is expressed by a natural number, then each of these distances is greater than 1.**

As in the previous theorem, if the distance between arbitrary two vertexes of a Diophantine *n*-polygon is equal to 1, then the its remaining vertexes should be located on the median of connecting these two vertexes segment. If we take into account that none of these vertexes are located on a one straight line, then we will have only the following cases (see Fig. 5).

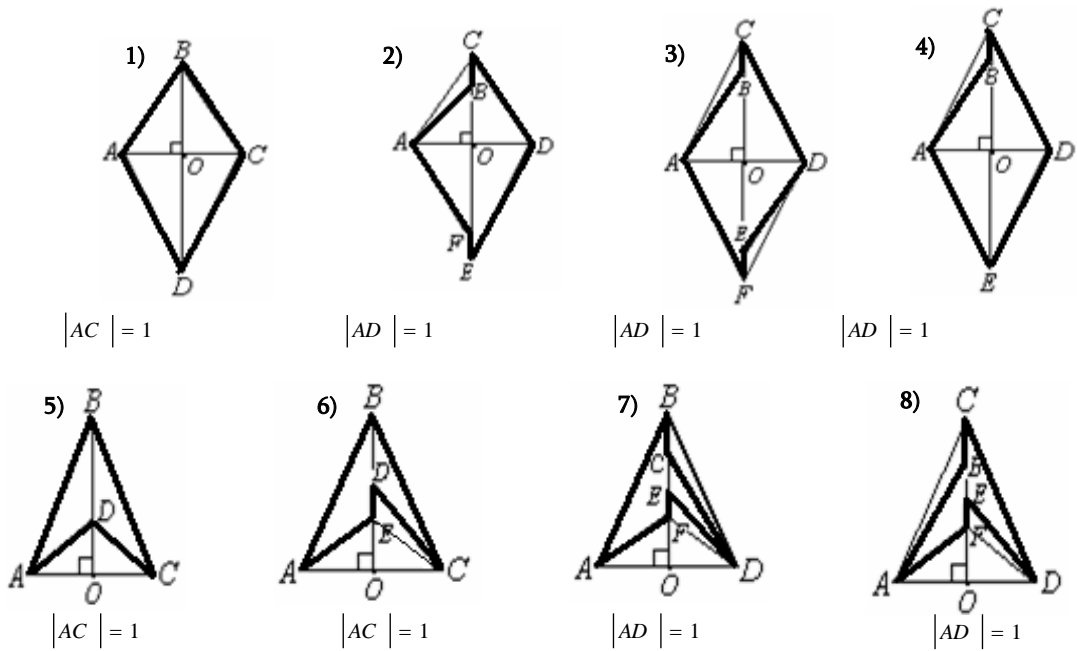

**Fig. 5**

The first fours have no solution according to **Theorem 1** all fours will be as 1) and since there is not exist such Diophantine rectangle, therefore are not exist 2), 3) and 4) Diophantine rectangles, and the last four have no solution by virtue of **Theorem 2** (Here, too 6), 7) and 8) by fill will be reduced to 5). I.e. our assumption is false or each of the distances given by the task conditions is greater than 1. Thereby fully is proved.

I.e. we have proved the following theorem.

From the tasks of studying bidiophantine geometric figures is obvious that one of the most important is the following.

**Task\*\* (n;k)** Exists or no for each fixed *k* integer the bidiophantine *n*-gon ($n > 3$), whose length of arbitrary side or diagonal is equal to *k*, and if exist, then find all such *n*.

At solution by us to the **task\* (n;1)**, it was shown that there is not exist such Diophantine n-gon (n> 3) those length of arbitrary side or diagonal is equal to 1. Proceeding from this, there can be not exist also such bidiophantine *n*-gon ($n>3$) whose length of arbitrary side or diagonal is equal to 1. As for the bidiophantine triangle, one of the sides of that is equal to 1, it should be isosceles accordingly of Lemma1, but then the vertex of this Diophantine triangle (the common point of the legs) cannot be Diophantine. I.e. there is not exist such bidiophantine triangle. This completely solved the task\*\* (*n*;1) for *k*=1, or

**Theorem\*\* (n;1) does not exist such a bidiophantine n-gon (n≥3), both convex and concave, those lengths of arbitrary sides or diagonals are equal to 1.**

Let us now show that even for k = 2 the **Task\*\* (*n;k*)** does not find a bidiophantine n-gon.

**Theorem\*\* (n;2)** does not find any bidiophantine *n*-gon ($n≥3$) whose distance between arbitrary two vertexes is equal to 2.

Let's assume the opposite. I.e. we assume that the length of each side and each diagonal of $A_1 A_2 A_3 ... A_n$ n-gon is expressed by a natural number, and at the same time the length of one of them is equal to 2. Then for all of triangles those base is this equal to 2 segment, and the third vertex of this n-polygon, the modulus of difference of leg's length is: 0 or 1 (this is easily illustrated by using the triangle inequality).

As there is not exist non-parallel to coordinate axis equal of 2 length bidiophantine segment (this leads due to the fact that there is no Pythagorean triangle with hypotenuse equal to 2). Therefore for each of the above mentioned triangles we have only one of the cases listed below (see Fig. 6 and Fig. 7).

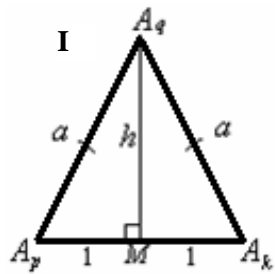
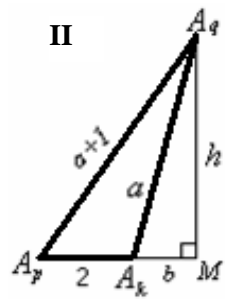

Fig. 6          Fig. 7

I. In this case, due that $[A_p A_k]$ and $\triangle A_p A_q A_k$ - are bidiophantines, therefore $\triangle A_p M A_q$ - also will be bidiophantine, but this is not possible, because there is not exist such Diophantine rectangular triangle those length of any cathetus is equal to 1, so if there is exist such bidiophantine triangle then it must not be isosceles.

II. In this case from the bidiophantine $\triangle A_p M A_q$ and $\triangle A_k M A_q$ due the Pythagorean theorem we have:

$$\begin{cases} a^2 = b^2 + h^2; \\ (a+1)^2 = (b+2)^2 + h^2; \\ a,b,h \in N. \end{cases} \Leftrightarrow \begin{cases} a^2 = b^2 + h^2; \\ a^2 + 2a + 1 = b^2 + 4b + 4 + h^2; \\ a,b,h \in N. \end{cases} \Leftrightarrow \begin{cases} a^2 = b^2 + h^2; \\ 2a + 1 = 4b + 4; \\ a,b,h \in N. \end{cases}$$

But this is impossible.

**So our assumption is false, or is not exist such bidiophantine *n*-gon (*n* > 3) both convex and concave, the length of any side or diagonal of that is equal to 2. Q.E.D.**

Now let's consider the **task\*\* (n;k)** for $k = 3$.

**Task\*\* (*n*;3)** Is there exist such bidiophantine n-**polygon** ($n > 3$) whose length of side or diagonal is equal to 3.

Let's assume that is exist such bidiophantine n-gon ($n > 3$) whose length of any side or diagonal is equal to 3. Then for all the triangles whose base represents this equal to 3 segment and the third vertex - any vertex of this n-gon, the modulus of difference of leg's length is equal to 0; 1 or 2 (this is easily illustrated by using the triangle inequality).

As there is not exist non-parallel to coordinate axis equal of 3 length bidiophantine segment, therefore for each of such triangles we have only one of the cases listed below (see Fig. 8, Fig. 9, Fig. 10).

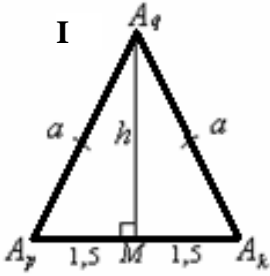 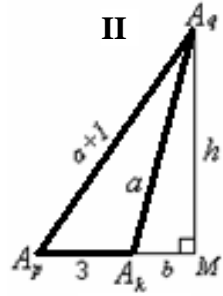 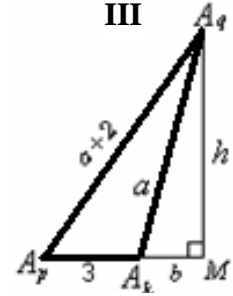

Fig. 8  Fig. 9  Fig. 10

I. In this case because $[A_p A_k]$ and $\triangle A_p A_q A_k$ - are bidiophantine, therefore $\triangle A_p M A_q$ will be bidiophantine. But this is impossible, because $[A_p M]$ is not bidiophantine.

II. In this case from the bidiophantine $\triangle A_p M A_q$ and $\triangle A_k M A_q$ due the Pythagorean theorem we have:

$$\begin{cases} a^2 = b^2 + h^2; \\ (a+1)^2 = (b+3)^2 + h^2; \\ a,b,h \in N. \end{cases} \Leftrightarrow \begin{cases} a^2 = b^2 + h^2; \\ 2a + 1 = 6b + 9; \\ a,b,h \in N. \end{cases} \Leftrightarrow \begin{cases} h^2 = 2(b+1)(b+2); \\ a = 3b + 4; \\ b,h \in N. \end{cases} \quad (1)$$

III. In this case from the bidiophantine $\triangle A_p MA_q$ and $\triangle A_k MA_q$ due the Pythagorean theorem we have:

$$\begin{cases} a^2 = b^2 + h^2; \\ (a+2)^2 = (b+3)^2 + h^2; \\ a,b,h \in N. \end{cases} \Leftrightarrow \begin{cases} a^2 = b^2 + h^2; \\ a^2 + 4a + 4 = b^2 + 6b + 9 + h^2; \\ a,b,h \in N. \end{cases} \Leftrightarrow \begin{cases} a^2 = b^2 + h^2; \\ 2(2a - 3b - 4) = 1; \\ a,b,h \in N. \end{cases} \quad (2)$$

But that is impossible.

Therefore for all such triangles $[A_p A_k]$ is parallel to the coordinate axis and

$$|A_p A_k| = 3; \; |A_p A_q| = 3b + 5; \; |A_q A_k| = 3b + 4; \; |A_k M| = b; \; |A_p M| = h = 2\sqrt{2(b+1)(b+2)},$$

where $b \in N$ and $2(b+1)(b+2)$ represents the square of the natural number.

$$\text{For } b = 2p \quad 2(b+1)(b+2) = 4(2p+1)(p+1) = t^2. \tag{3}$$

$$\text{and for } b = 2p - 1 \text{ the } 2(b+1)(b+2) = 4p(2p+1) = t^2 \tag{4}$$

(3) and (4) have infinitely many solutions in natural numbers. Indeed, because $(p+1; 2p+1) = 1$ therefore (3) will be found such $m, n \in N$, for that

$$\begin{cases} p + 1 = m^2; \\ 2p + 1 = n^2. \end{cases} \Rightarrow n^2 - 2m^2 = -1. \tag{5}$$

$$\begin{cases} (a_1^2 - 2b_1^2)(c_1^2 - 2d_1^2) = (a_1 c_1 + 2b_1 d_1)^2 - 2(a_1 d_1 + b_1 c_1)^2 & (6) \\ (a_1^2 - 2b_1^2)(c_1^2 - 2d_1^2) = (a_1 c_1 - 2b_1 d_1)^2 - 2(a_1 d_1 - b_1 c_1)^2 & (7) \end{cases}$$

Due the (6) and (7) by combining solutions of equations $c_1^2 - 2d_1^2 = 1$ and $a_1^2 - 2b_1^2 = -1$ we can obtain all new solutions of (5). Below are presented some results of computations:

| $n$ | $m$ | $n^2 - 2m^2$ | $b$ | $2\sqrt{2(b+1)(b+2)}$ |
|---|---|---|---|---|
| 7 | 5 | $-1$ | 48 | 140 |
| 41 | 29 | $-1$ | 1680 | 4756 |
| 239 | 169 | $-1$ | 57120 | 80782 |
| 1393 | 985 | $-1$ | 1 940 448 | 2 744 210 |
| 8119 | 5741 | $-1$ | 65 918 160 | 93 222 358 |

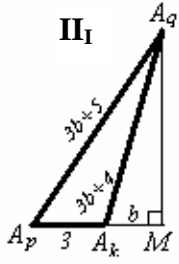

ნახ. 11

It is much easier to obtain the solutions of (4). Indeed, as $(p\,;\,2p+1)=1$, therefore, we can more easily obtain the solution of (4). It follows that there will be existing such $m_1, n_1 \in N$ for that

$$\begin{cases} p = m_1^2; \\ 2p+1 = n_1^2. \end{cases} \Rightarrow n_1^2 - 2m_1^2 = 1. \qquad (8)$$

Thus we can obtain the solution of (4) by solving Pell's equation (8). In addition $b = 2m_1^2 - 1$.

II$_I$. type set of biophysical triangles, as it is shown above, is not an empty set.

$$\cos \widehat{A_p A_K A_q} = -\cos \widehat{A_q A_K M} = -\frac{b}{3b+4}. \qquad \cos \widehat{A_q A_p A_k} = \frac{b+3}{3b+5}.$$

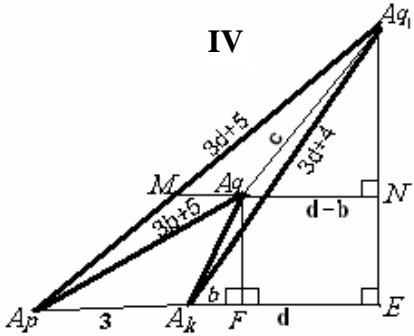

ნახ. 12

Therefore by increasing $b$ decreases $\cos \widehat{A_p A_k A_q}$ and $\cos \angle \widehat{A_q A_p A_k}$. I.e. are increasing $\widehat{A_p A_k A_{q_1}}$ and $\widehat{A_q A_p A_{k_1}}$.

Therefore, if $b, d \in N$ $(d > b)$, then, for type II$_I$ triangles we would have case IV or V.

Let's consider case IV.

From $\Delta A_p A_q A_{q_1}$ $c > (3d+5) - (3b+5) = 3(d-b)$

$$A_q \in [MN] \| [A_p A_k] \Rightarrow \Delta MNA_q \sim \Delta A_p EA_{q_1} \Rightarrow \frac{|MA_{q_1}|}{|A_p A_{q_1}|} = \frac{|NA_{q_1}|}{|EA_{q_1}|} \Rightarrow$$

$$\Rightarrow |MA_{q_1}| = (3d+5)\left(\frac{2\sqrt{2(d+1)(d+2)} - 2\sqrt{2(b+1)(b+2)}}{2\sqrt{2(d+1)(d+2)}}\right) =$$

$$= (3d+5)\left(1 - \sqrt{\frac{b+1}{d+1} \cdot \frac{b+2}{d+2}}\right) < (3d+5)\left(1 - \frac{b+1}{d+1}\right) = \frac{(3d+5)(d-b)}{d+1} =$$

$$= 3(d-b) + \frac{2(d-b)}{d+1} < 3(d-b) + 2.$$

$$\begin{cases} 3(d-b) < c < 3(d-b)+2; \\ b, d \in N. \end{cases} \Rightarrow c = 3(d-b)+1.$$

$\Delta A_q NA_{q_1}$ -დან $|A_q A_{q_1}|^2 = |A_q N|^2 + |NA_{q_1}|^2 \Rightarrow$

$$\Rightarrow (3(a-b)+1)^2 = (d-b)^2 + 8\left(\sqrt{(d+1)(d+2)} - \sqrt{(b+1)(b+2)}\right)^2 \Leftrightarrow$$

$$\Leftrightarrow 8(d-b)^2 + 6(d-b) + 1 = 8\left(\sqrt{(d+1)(d+2)} - \sqrt{(b+1)(b+2)}\right)^2,$$

that is impossible, because the left side of the last equation is odd, and the right side if it is natural, ☐ will be even.

Now let's consider the case V. It is obvious that $|A_qO| \in (3b+4; 3b+5)$ and $|A_lO| \in (3d+4; 3d+5)$, therefore $|A_qA_l| \in (3b+3d+8; 3b+3d+10)$ and if $|A_qA_l| \in N$, then $|A_qA_l| = 3b+3d+9 = 3(b+d+3)$.

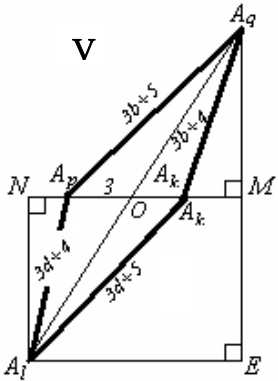
Fig. 13

(From the rectangle $A_lNME$) $\begin{cases} |A_lE| = |NM| = b+d+3; \\ |ME| = |A_lN|. \end{cases}$,

$(|A_qM|, |A_lN| \in N) \Rightarrow |A_qE| \in N$, as $|A_qE| = |A_qM| + |A_lN|$.

From the rectangle triangle $\triangle A_lEA_q$ $|A_lA_q|^2 = |A_lE|^2 + |A_qE|^2$, from that we will obtain $9(b+d+3)^2 = (b+d+3)^2 + |A_qE|^2 \Leftrightarrow 8(b+d+3)^2 = |A_qE|^2$. This is impossible, because $b, d, |A_qE| \in N$. I.e. are not existing IV and V type bidiophantine rectangles, so all bidiophantine *n*-gons those length of side is equal to 3 should only be of type VI (on one side of a straight line containing equal to 3 the side length) and whose diagonal length is equal to 3 would be only of type VII only.

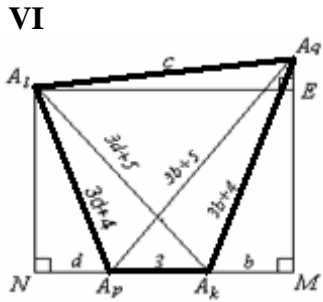

Fig. 14

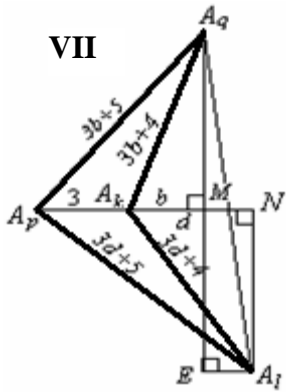

Fig. 15

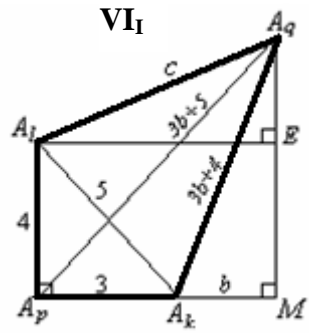

Fig. 16

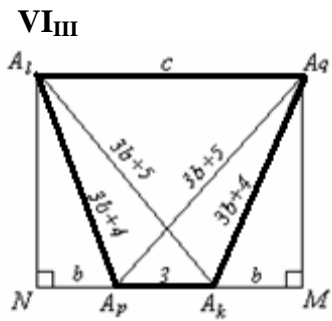

Fig. 17

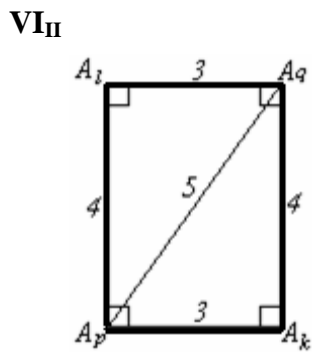

Fig. 18

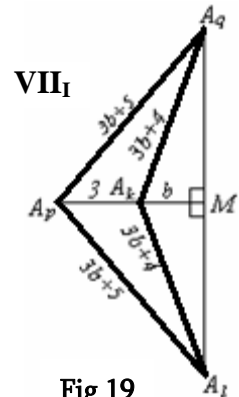

Fig 19

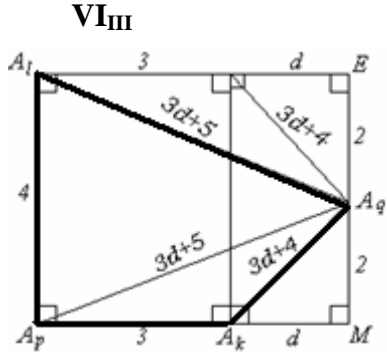
**VI<sub>III</sub>**

Fig. 20

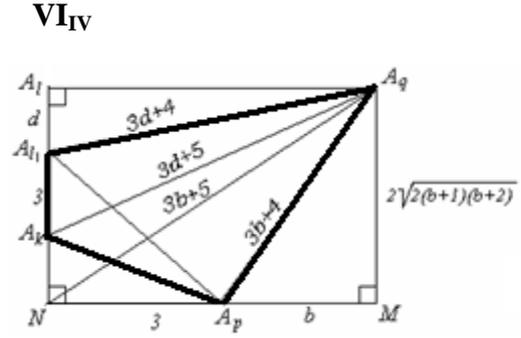
**VI<sub>IV</sub>**

Fig. 21

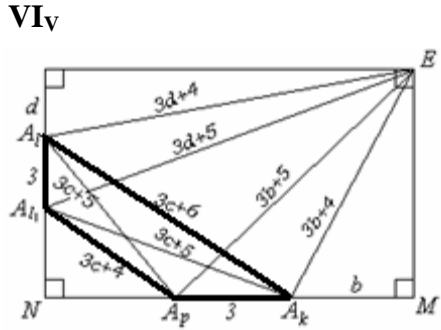
**VI<sub>V</sub>**

Fig. 22

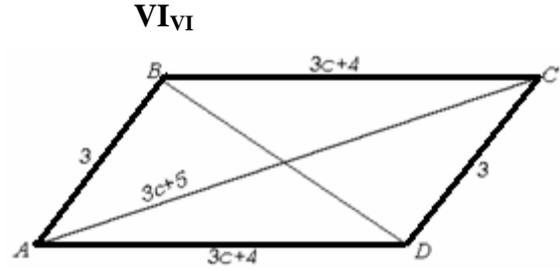
**VI<sub>VI</sub>**

Fig. 23

For the case VI from $\triangle A_l E A_q$ due the Pythagorean Theorem: we have: $|A_l A_q|^2 = |A_l E|^2 + |A_q E|^2$.

from that we obtain

$$c^2 = (b + d + 3)^2 + 4\left(\sqrt{2(b+1)(b+2)} - \sqrt{2(d+1)(d+2)}\right)^2. \quad (9)$$

Let's consider (9) for $d = 0$ (this is the case VI$_1$). For this case we have:

$|A_P A_l| = 4$, $|A_P A_K| = 3$, $|A_l A_K| = 5$, $|A_P A_q| = 3b + 5$, $|A_K A_q| = 3b + 4$, $|A_K M| = b$,

$|A_q M| = 2\sqrt{2(b+1)(b+2)}$, $|A_q E| = 2\sqrt{2(b+1)(b+2)} - 4$. Therefore teh (9) will be as:

$$c^2 = (b + 3)^2 + 4\left(\sqrt{2(b+1)(b+2)} - 2\right)^2. \quad (10)$$

From the blunt triangle $\triangle A_P A_l A_q$ due the inequality of triangles we have:

$$\begin{cases} c + 4 > 3b + 5; \\ c < 3b + 5. \end{cases} \Rightarrow c \in (3b + 1; 3b + 5).$$

As $c \in N$, thus $c \in \{3b + 2; 3b + 3; 3b + 4\}$. (11)

By transformation of (10) we will have:

$$16\sqrt{2b^2 + 6b + 4} = 9b^2 + 30b + 41 - c^2. \quad (12)$$

Due the taking into account (11) and (12) we will have:

$$\begin{cases} \left[\begin{array}{l} \left[16\sqrt{2b^2+6b+4} = 9b^2+30b+41-(3b+2)^2;\right. \\ \left. 16\sqrt{2b^2+6b+4} = 9b^2+30b+41-(3b+3)^2;\right. \\ \left. 16\sqrt{2b^2+6b+4} = 9b^2+30b+41-(3b+4)^2.\right. \end{array}\right. \\ b \in Z_0. \end{cases} \Leftrightarrow \begin{cases} \left[\begin{array}{l} 16\sqrt{2b^2+6b+4} = 18b+37; \\ 16\sqrt{2b^2+6b+4} = 12b+32; \\ 16\sqrt{2b^2+6b+4} = 6b+25. \end{array}\right. \\ b \in Z_0. \end{cases} \quad (13)$$

The first and third equations of the system (13) did not have the solution in positive integers, because the left side of these equations, if it is natural, is even, as well as the right side is odd. Thus from (13) we will obtain:

$$\begin{cases} 16\sqrt{2b^2+6b+4} = 12b+32; \\ b \in Z_0. \end{cases} \Leftrightarrow \begin{cases} 4\sqrt{2b^2+6b+4} = 3b+8; \\ b \in Z_0. \end{cases} \Leftrightarrow$$

$$\Leftrightarrow \begin{cases} 16(2b^2+6b+4) = 9b^2+48b+64; \\ b \in Z_0. \end{cases} \Leftrightarrow \begin{cases} 23b^2+48b = 0; \\ b \in Z_0. \end{cases} \Leftrightarrow \begin{cases} \left[\begin{array}{l} b = -\dfrac{48}{23}; \\ b = 0. \end{array}\right. \\ b \in Z_0. \end{cases} \Leftrightarrow (b=0).$$

In this case we obtain a VI$_{II}$ type rectangle.

Now let's consider (9) for $d = b \neq 0$. In this case we have VI$_{III}$. Such bidiophantine rectangles exist for those $b \in N$, for that $\sqrt{2(b+1)(b+2)}$ - is a natural number. As it is mentioned above, such numbers are numerous.

As for (9), its solution in case when $b, d \in N$ and $b \neq d$ by us in **Basic** software was tested from $b$ and $d$ up to 1,000,000, but were not find such natural $c$, $b$ and $d$.

In our opinion such $c$, $b$ and $d$ do not exist. **It is noteworthy that this is possible to be considered as a task on a bidiophantine ellipse, those focuses are located in $A_l$ and $A_q$ points and length of major diameter is equal to $3b+3d+9$. Therefore, the existence of such bidiophantine rectangle is equivalent to the existence on a bidiophantine ellipse of a bidiophantine cord that is not parallel to the coordinate axis.**

Similarly to (9) for VII from the $\triangle A_q E A_l$ rectangle we obtain:

$$c^2 = (d-b)^2 + 4\left(\sqrt{2(b+1)(b+2)} + \sqrt{2(d+1)(d+2)}\right)^2. \quad (14)$$

For $b = d$ we obtain case VIII and (14) will be as

$$\begin{cases} c^2 = 4\left(\sqrt{2(b+1)(b+2)}\right)^2; \\ c, b \in N. \end{cases} \Leftrightarrow \begin{cases} c = 2\sqrt{2(b+1)(b+2)}; \\ c, b \in N. \end{cases} \quad (15)$$

As it is shown above the (15) has lots of solutions, and as for (14), similarly to (9) has been tested on **Basic** software was tested from $b$ and $d$ up to 1,000,000, but were not find such natural $c$, *b* and *d* and up to 1,000,000, but no such c, and not found. (In our opinion there is no such rectangle).

VI$_{III}$ has not solution because the value of [$EA_q$] will be natural.

VI$_{IV}$ has not solution because $\triangle A_{l_1}NA_p$ and $\triangle A_k NA_p$ at the same time will be Diophantine right triangles that is impossible because the Diophantine equation $x^2+3^2=y^2$ has the unique solution in natural numbers $x=4$ and $y=5$.

VI$_V$ is impossible because $|A_{l_1}A_k|+|A_l A_p|=|A_{l_1}A_p|+|A_l A_k|$, this contradicts the well-known task accordingly of that sum of the lengths of diagonals of convex rectangle will be greater than the sum of the lengths of the opposite sides.

VI$_{VI}$ is impossible because we would then have that $|BD|=3c+3$ or $|BD|=3c+5$.

The first case ($|BD|=3c+3$) is impossible because we would have that $|BD|+|AC|=|BC|+|AD|$ that contradicts to the mentioned in VI$_V$ provision.

In the second case we obtain that ΔBAD is a right triangle and $c = 0$ (this was considered in VI$_{II}$).

Let's get back to V.

$$\begin{cases} \cos A_q \widehat{A}_p A_k = \dfrac{b+3}{3b+5} = \dfrac{1}{3}+\dfrac{4}{3(3b+5)}; \\ \cos A_l \widehat{A}_p A_k = -\dfrac{d}{3d+4}=-\dfrac{1}{3}+\dfrac{4}{3(3d+4)}; \\ d,b \in N. \end{cases} \Rightarrow \begin{cases} A_q \widehat{A}_p A_k \in (0;\arccos \dfrac{1}{3}); \\ A_l \widehat{A}_p A_k \in (0;\pi - \arccos \dfrac{1}{3}) \end{cases} \Rightarrow A_l \widehat{A}_p A_q \in (0;\pi).$$

Similarly $A_l \widehat{A}_k A_q \in (0;\pi)$.

**We obtain that for *k*=3, the bidiophantine n-gon will only be either a triangle, or a rectangle, because it cannot have the following shapes (see Fig. 24, Fig. 25, Fig. 26).**

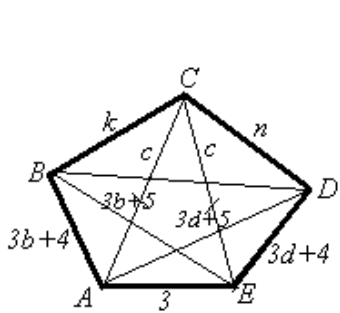 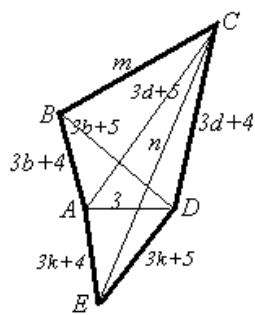 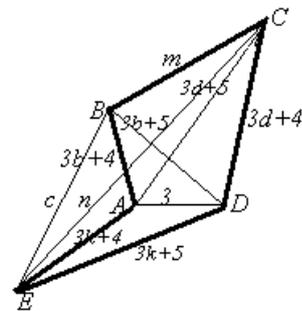

**Fig. 24**      **Fig. 25**      **Fig. 26**

Finally we have

**Theorem\*\* (n;3)** If for bidiophantine *n*-gon $k = 3$, then $n \in \{3; 4\}$.

Now let's consider the **Task\*\* (n;k)** for $k = 4$.

**Task\*\* (n;4) Is there existing or not such bidiophantine n-gon ($n \geq 3$) those side or diagonal length is equal to 4, and if exists let's find all such $n$.**

Let's suppose that there is existing such bidiophantine $n$-gon ($n \geq 3$) whose length of any side or diagonal is equal to 4. Then for all those triangles the base of that is represented by this equal to 4 length of segment, and the vertex is any of these $n$-polygon's vertex the modulus of difference of a legs lengths is equal to 0; 1; 2 or 3 (This is easily to show by using the triangle inequality).

As there is not existing non-parallel to the coordinate axis bidiophantine segment those length is equal to 4, (this implies that there is no Pythagorean triangle with hypotenuse length equal to 4), we have one of the following four cases for all such triangles (see Fig. 27, Fig. 28, Fig. 29, Fig. 30).

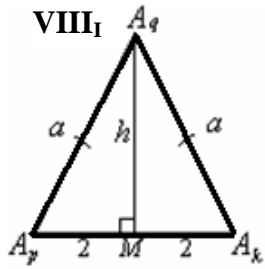
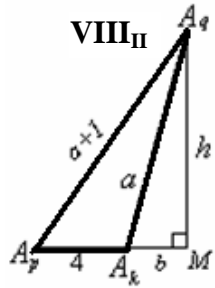
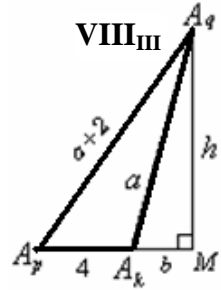
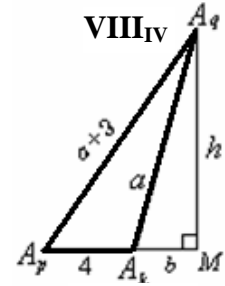

Fig. 27    Fig. 28    Fig. 29    Fig. 30

In the case VIII$_I$, as $[A_p A_k]$ and $\Delta A_p A_q A_k$ - are bidiophantine, thus $\Delta A_p M A_q$ also will be bidiophantine, but this is impossible because is not exist bidiophantine right triangle, the length of those cathetus is 2.

In the case VIII$_{II}$ we have from the bidiophantine triangles $\Delta A_p M A_q$ and $\Delta A_k M A_q$ accordingly of Pythagorean theorem:

$$\begin{cases} a^2 = b^2 + h^2; \\ (a+1)^2 = (b+4)^2 + h^2; \\ a,b,h \in N. \end{cases} \Leftrightarrow \begin{cases} a^2 = b^2 + h^2; \\ a^2 + 2a + 1 = (b^2 + h^2) + 8b + 16; \\ a,b,h \in N. \end{cases} \Leftrightarrow \begin{cases} a^2 = b^2 + h^2; \\ 2a + 1 = 8b + 16; \\ a,b,h \in N. \end{cases}$$

This is impossible.

In the case VIII$_{IV}$, similarly to VIII$_{II}$, from the right bidiophantine triangles $\Delta A_p M A_q$ and $\Delta A_k M A_q$ we have: $\begin{cases} a^2 = b^2 + h^2; \\ 6a + 9 = 8b + 16; \\ a,b,h \in N. \end{cases}$ this is impossible.

In the case VIII$_{III}$ bidiophantine triangles $\triangle A_p M A_q$ and $\triangle A_k M A_q$ accordingly of Pythagorean theorem we have:

$$\begin{cases} a^2 = b^2 + h^2; \\ (a+2)^2 = (b+4)^2 + h^2; \\ a,b,h \in N. \end{cases} \Leftrightarrow \begin{cases} a^2 = b^2 + h^2; \\ a = 2b + 3; \\ a,b,h \in N. \end{cases} \Leftrightarrow \begin{cases} h^2 = 3(b+1)(b+3); \\ a = 2b + 3; \\ b,h \in N. \end{cases}$$

So for all such triangles $[A_p A_K]$ are parallel to the coordinate axis. In addition: $|A_p A_K| = 4$; $|A_p A_q| = 2b + 5$; $|A_q A_K| = 2b + 3$; $|A_q M| = h = \sqrt{3(b+1)(b+3)}$, where $b \in N$ and $3(b+1)(b+3)$ represents the full square of the natural number. We tested up to 25,000 and such $b$ are found 54 in total.

| $b$ | $\sqrt{3(b+1)(b+3)}$ | $b$ | $\sqrt{3(b+1)(b+3)}$ |
|---|---|---|---|
| 5 | 12 | 5820 | 10084 |
| 24 | 45 | 7171 | 12424 |
| 95 | 168 | 7951 | 13775 |
| 360 | 627 | 8731 | 15126 |
| 1349 | 2340 | 9511 | 16477 |
| 5040 | 8733 | 10082 | 17466 |
| 5820 | 10084 | 10862 | 18817 |
| 7171 | 12424 | 11433 | 19806 |
| 7951 | 13775 | 11642 | 20168 |

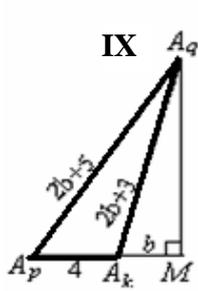
IX

Fig. 31

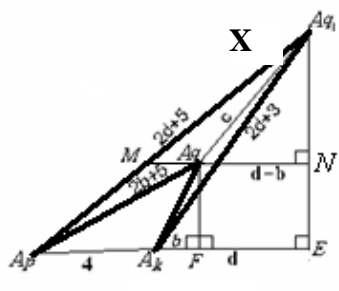
X

Fig. 32

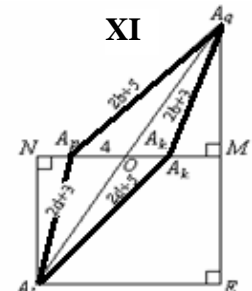
XI

Fig. 33

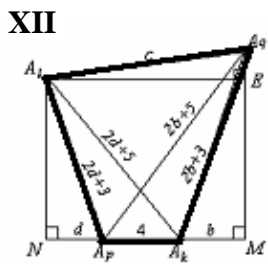
XII

Fig. 34

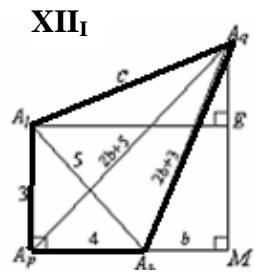
XII$_I$

Fig. 35

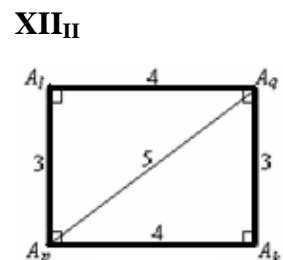
XII$_{II}$

Fig. 36

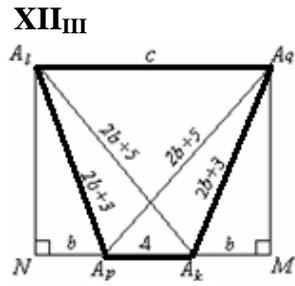
Fig. 37

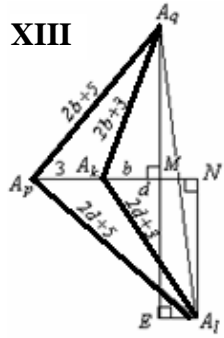
Fig. 38

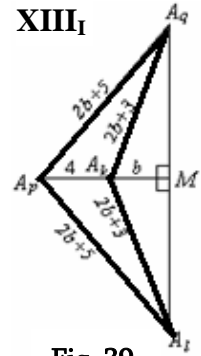
Fig. 39

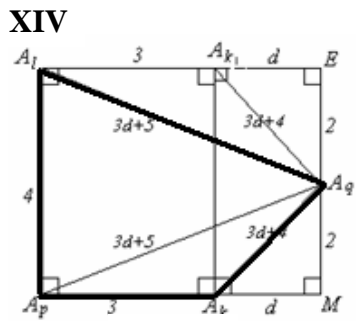
Fig. 40

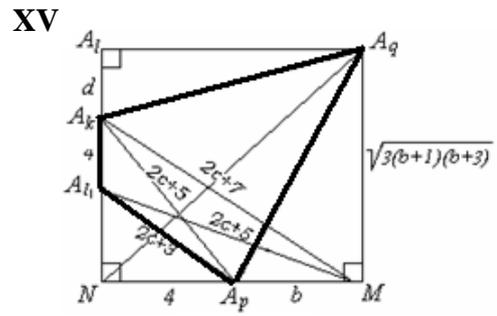
Fig. 41

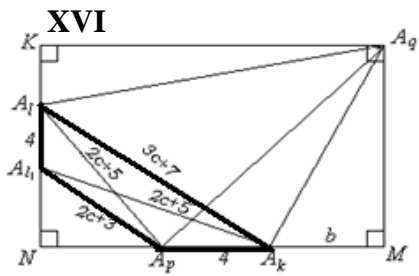
Fig. 42

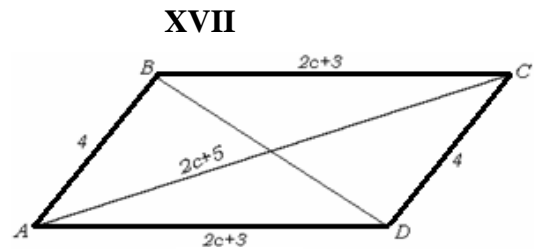
Fig. 43

The set of type IX Bidiophantine triangles, as it is shown above, is not an empty set.

$$\begin{cases} \cos \widehat{A_q A_p A_k} = \dfrac{b+4}{2b+5}; \\ \cos \widehat{A_p A_K A_q} = -\cos \widehat{A_q A_K M} = -\dfrac{b}{2b+3} \end{cases}.$$

Therefore, by increasing in $b$ decreases $\cos \widehat{A_p A_K A_q} = -\cos \widehat{A_q A_K M} = -\dfrac{b}{2b+3}$.

Here, similarly to the previous paragraph, by the increasing of $b$ decreases $\cos \widehat{A_p A_K A_q}$ and $\cos A_q \widehat{A}_p A_k$. I.e. increases $A_p \widehat{A}_k A_q$ and $A_q \widehat{A}_p A_k$. Therefore for $b, d \in N$ $(d > b)$ type IX triangles we can only have cases X and XI.

Let's consider the case X

From $\Delta A_p A_q A_{q_1}$ $c > (2d + 5) - (2b + 5) = 2(d - b)$.

$[MN] \| [A_p A_k] \Rightarrow \Delta MNAq_1$  $\Delta A_p EA\varnothing_\perp \Rightarrow \dfrac{|MA_{q_1}|}{2d+5} = \dfrac{|NA_{q_1}|}{|EA_{q_1}|} \Rightarrow |MA_{q_1}| =$

$= (2d + 5)\left(\dfrac{\sqrt{3(d+1)(d+3)} - \sqrt{3(b+1)(b+3)}}{\sqrt{3(d+1)(d+2)}}\right) = (2d+5)\left(1 - \sqrt{\dfrac{b+1}{d+1} \cdot \dfrac{b+3}{d+3}}\right) <$

$< (2d+5)\left(1 - \dfrac{b+1.5}{d+1.5}\right) = \dfrac{(2d+5)(d-b)}{d+1.5} = 2(d-b) + \dfrac{d-b}{d+1.5} < 2(d-b) + 2;$

$\begin{cases} 2(d-b) < c < 2(d-b) + 2; \\ b, d, c \in N. \end{cases} \Rightarrow c = 2(d-b) + 1.$

From the right triangle $\Delta A_q NA_{q_1}$ $|A_q A_{q_1}|^2 = |A_q N|^2 + |NA_{q_1}|^2$ i.e.

$(2(d-b) + 1)^2 = (d-b)^2 + 3\left(\sqrt{(d+1)(d+3)} - \sqrt{(b+1)(b+3)}\right)^2 \Leftrightarrow$

$\Leftrightarrow 2(3db + 4d + 8b + 8) + 1 = 6\sqrt{(b+1)(b+3)(d+1)(d+3)}.$  (16)

This is impossible, because its left side is odd, and its right side if it will be natural - even.

From the case XI is obvious that

$\begin{cases} |A_q O| \in (2b+3; 2b+5); \\ |A_l O| \in (2d+3; 2d+5). \end{cases} \Rightarrow |A_l A_q| \in (2b+2d+6; 2b+2d+10).$

In addition, if we take into account that from $\Delta A_l A_q A_K$ $|A_l A_q| < |A_l A_K| + |A_q A_K| =$

$= 2d + 5 + 2b + 3 = 2b + 2d + 8$, and $|A_q A_l| \in N$, then we obtain $|A_q A_l| = 2b + 2d + 7$.

From the right triangle $\Delta A_l EA_q$ accordingly of Pythagorean theorem we have $|A_l A_q|^2 = |A_l E|^2 + |A_q E|^2$ from that we derive

$(2b + 2d + 7)^2 = (b + d + 4)^2 + 3\left(\sqrt{(b+1)(b+3)} + \sqrt{(d+1)(d+3)}\right)^2 \Leftrightarrow$

$\Leftrightarrow 6bd + 8(b+d) + 15 = 6\sqrt{(b+1)(b+3)(d+1)(d+3)}.$  (16$_1$).

(16$_1$) has no solution in natural numbers, because the left side is odd and the right is even

In the case XII from the right triangle $\Delta A_l EA_q$ accordingly of Pythagorean theorem we have $|A_l A_q|^2 = |A_l E|^2 + |A_q E|^2$ from that we derive

$$c^2 = (b + d + 4)^2 + 3\left(\sqrt{(b + 1)(b + 3)} + \sqrt{(d + 1)(d + 3)}\right)^2. \tag{17}$$

In (17) for $d=0$ (XIII case) we will have:

$$|A_p A_l| = 3; |A_p A_k| = 4; |A_l A_k| = 5; |A_p A_q| = 2b + 5; |A_k A_q| = 2b + 3;$$

$$|A_k M| = b; |A_q M| = \sqrt{3(b + 1)(b + 3)}; |A_q E| = \sqrt{3(b + 1)(b + 3)} - 3.$$

and in this case (17) will be as:

$$c^2 = (b + 4)^2 + \left(\sqrt{3(b + 1)(b + 3)} - 3\right)^2. \tag{18}$$

From the blunt triangle $\Delta A_p A_l A_q$, due the inequality of triangles we have:

$$\begin{cases} c + 3 > 2b + 5; \\ c < 2b + 5. \end{cases} \Rightarrow c \in (2b + 2; 2b + 5), \text{ but as } c \in N, \text{ thus } c \in \{2b + 3; 2b + 4\}. \tag{19}$$

By transformation of (18) we obtain $6\sqrt{3b^2 + 12b + 9} = 4b^2 + 20b + 34 - c^2$. (20)

The left side of (18) must necessarily be a natural number and then it will be even, thus in order for the right side also to be even, from (19) would be used only $c = 2b + 4$ (21) i.e. by substitution of (21) in (20) we will obtain:

$$\begin{cases} 6\sqrt{3b^2 + 12b + 9} = 4b^2 + 20b + 34 - (2b + 4)^2; \\ b \in Z_0. \end{cases} \Rightarrow \begin{cases} 3b^2 + 144b = 0; \\ b \in Z_0. \end{cases} \Leftrightarrow b = 0.$$

In this case we obtain the $XII_{II}$ rectangle.

Now let's consider (17) for $d = b \neq 0$.

In this case we have $XII_{III}$. Such bidiophantine rectangles exist only such $b \in N$ for that $\sqrt{3(b + 1)(b + 3)}$ - is a natural number. As we have shown, such numbers are numerous.

As for the solutions of (17), in the case when $b, d \in N$ and $b \neq d$, by us in **Basic** software was tested from $b$ and $d$ up to 1,000,000, but were not find such natural $c$, $b$ and $d$. In our opinion such $b$ and $d$ do not exist. **It is noteworthy that this case is possible to be considered as a task on bidiophantine ellipse, those focuses are located in $A_l$ and $A_q$ points and those major diameter is equal to 2b+2d+8. Therefore, the existence of such $A_p A_l A_q A_k$ bidiophantine rectangle is equal to the existence on an ellipse of such bidiophantine cord that is not parallel to the coordinate axis.**

Similarly of (17) in case XIII from right triangle $\Delta A_q E A_l$ we have

$$c^2 = (d - b)^2 + 3\left(\sqrt{(b + 1)(b + 3)} + \sqrt{(d + 1)(d + 3)}\right)^2. \tag{21}$$

For $b = d$ we obtain $XIII_I$ and (21) will be as

$$\begin{cases} c^2 = 4\left(\sqrt{3(b+1)(b+3)}\right)^2; \\ c, b \in N. \end{cases} \Leftrightarrow \begin{cases} c = 2\sqrt{3(b+1)(b+3)}; \\ c, b \in N. \end{cases} \quad (22)$$

As it was shown above, (22) has several solutions, and as for (21) for $b \neq d$ by us in **Basic** software was tested from $b$ and $d$ up to 1,000,000, but were not find such natural $c$, $b$ and $d$.

XIV does not have solution in natural numbers, because then $\Delta A_k MA_q$ must be bidiophantine, but there is not exist Pythagorean triangle with one of the cathetus equal to 2.

XV does not have solution in natural numbers because otherwise for the $\square MA_p A_{l_1} A_k$ would be violated the inequality

$$|A_k A_p| + |A_{l_1} M| > |A_{l_1} A_p| + |A_k M|.$$

The XVI case is impossible for the same reason as for XV.

In the XVII variant we have that $|BD|=2c+1$ or $|BD|=2c+5$. In the first case, the true inequality $|BD|+|AC|>|BC|+|AD|$ is violated that is impossible, and in the second case we obtain that $\Delta BAD$ is a right triangle and $c=0$ (this was considered in XIIII).

Let's go back to XI:

$$\begin{cases} \cos \widehat{A_q A_p A_k} = \dfrac{b+4}{2b+5} = \dfrac{1}{2} + \dfrac{1.5}{2b+5}; \\ \cos \widehat{A_l A_p A_k} = -\dfrac{d}{2d+3} = -\dfrac{1}{2} + \dfrac{1.5}{2d+3}. \end{cases} \Rightarrow \begin{cases} \widehat{A_q A_p A_k} \in \left(0; \dfrac{\pi}{3}\right); \\ \widehat{A_l A_p A_k} \in \left(0; \dfrac{2\pi}{3}\right). \end{cases} \Rightarrow \widehat{A_l A_p A_q} \in (0; \pi)$$

Similarly $\widehat{A_l A_k A_q} \in (0; \pi)$.

We obtain that for $k=4$ the bidiophantine $n$-polygon would be either a only triangle or a rectangle because it cannot have the following shapes:

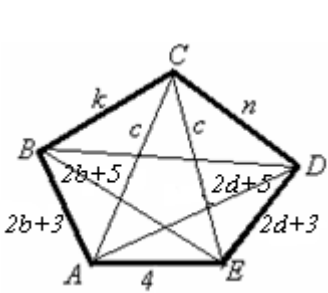 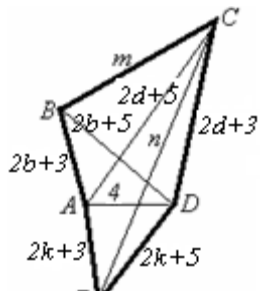 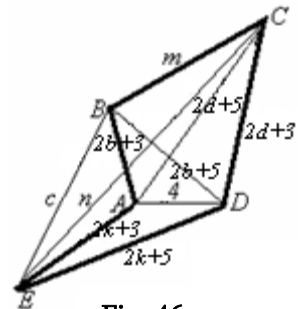

**Fig. 44**     **Fig. 45**     **Fig. 46**

Finally we have:

**Theorem** (n;4) **If for bidiophantine $n$-gon $k = 4$, then $n \in \{3; 4\}$.**

We investigate the **Task** (n;k) as well as for $k = 5$, $k = 6$, $k = 7$. We have some cases to checked. The reader will find them in the following work.

To the question whether there exist for each $k \in N$ $(k \geq 3)$ bidiophantine triangle, the length of any side of that is equal to *k*, the answer is positive. See Fig. 47.

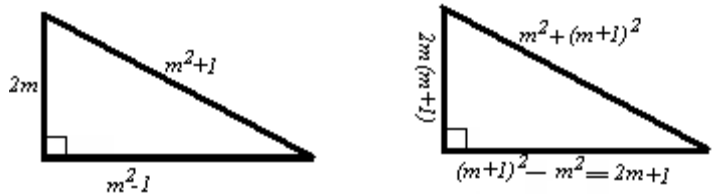

**Fig. 47**

To the question whether there exist a bidiophantine rectangle for each of $k \in N$ $(k \geq 3)$ the length of any sides of those is *k*, the answer is also positive. See Fig. 48

$$m \in N \setminus \{1\}$$

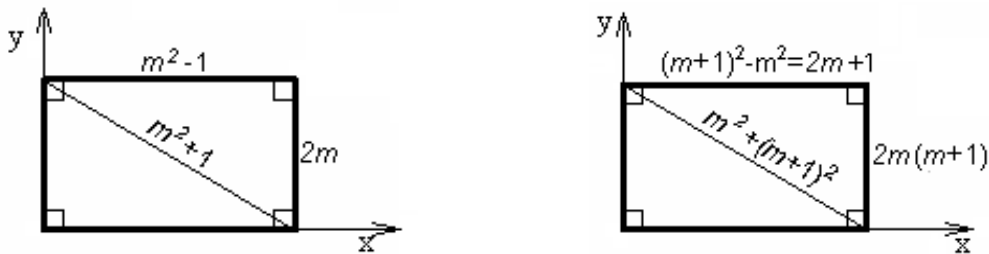

Fig.48

CONCLUSIONS

The goal of the work is to take on and study one of the fundamental tasks studying Bidiophantine n-gons (the author of the paper considers an integral n-gon is Bidiophantine if the coordinates of its each vertex are integer numbers).

Task **(*n; k*): is there a Bidiophantine n-gon ($n \geq 3$) with any side or diagonal equal to *k* for each fixed natural number *k* ? In case it exists then let us find each such *n*.

As a result of fundamental studies we have obtained a full answer to the above mentioned question for k=1, k=2, k=3 and k=4.

The paper shows that for k=1 and k=2 such Bidiophantine n-gon does not exist, and for k=3 and k=4 definitely n is equal to 3 or 4. It is also shown that when $k>2$ there is always a Bidiophantine rectangle with the length of one of the sides equal to k